\documentclass[11pt,reqno]{article}

\usepackage[usenames]{color}
\usepackage{amssymb}
\usepackage{graphicx}
\usepackage{amscd}
\usepackage[english]{babel}
\usepackage{color}
\usepackage{fullpage}
\usepackage{psfig}
\usepackage{graphics,amsmath,amssymb}
\usepackage{amsthm}
\usepackage{amsfonts}
\usepackage{latexsym}
\usepackage{epsf}
\usepackage{float}

\usepackage[colorlinks=true,
linkcolor=webgreen,
filecolor=webbrown,
citecolor=webgreen]{hyperref}

\definecolor{webgreen}{rgb}{0,.5,0}
\definecolor{webbrown}{rgb}{.6,0,0}

\newcommand{\beql}[1]{\begin{equation}\label{#1}}
\newcommand{\eeq}{\end{equation}}
\newcommand{\eqn}[1]{(\ref{#1})}

\newcommand{\Snj}{\scriptstyle}

\newtheorem{thm}{Theorem}{\bfseries}{\itshape}
{\bfseries}{\itshape}
\newtheorem{lem}[thm]{Lemma}{\bfseries}{\itshape}
{\bfseries}{\itshape}
\newtheorem{conj}{Conjecture}{\bfseries}{\itshape}

\newcommand{\LIMA}{{47}}
\newcommand{\LIMB}{{48}} 
\newcommand{\LIMC}{{80}}

\begin{document}
\theoremstyle{plain}


\begin{center}
{\large\bf The Curling Number Conjecture\footnote{For a sequel to this paper, see arXiv:1212.6102}} \\
\vspace*{+.2in}

Benjamin Chaffin, \\
Intel Processor Architecture, \\
2111 NE 25th Avenue, Hillsboro, OR 97124, USA, \\
Email: \href{mailto:chaffin@gmail.com}{\tt chaffin@gmail.com} \\

and \\

N. J. A. Sloane,\footnote{To whom correspondence should be addressed} \\
The OEIS Foundation Inc., \\
11 South Adelaide Ave., Highland Park, NJ 08904, USA \\
Email: \href{mailto:njasloane@gmail.com}{\tt njasloane@gmail.com} \\


\vspace*{+.1in}
Dec. 11, 2009; revised Feb. 13, 2010, Sep. 24, 2012
\vspace*{+.1in}


{\bf Abstract}
\end{center}

Given a finite nonempty sequence of integers $S$,
by grouping adjacent terms it is always possible
to write it, possibly in many ways,
as $S = X\,Y^k$, where $X$ and $Y$ are sequences
and $Y$ is nonempty.
Choose the version which maximizes the value of $k$: this
$k$ is the {\em curling number} of $S$.
The {\em Curling Number Conjecture} is that if
one starts with {\em any} initial sequence $S$,
and extends it by repeatedly appending the curling number of 
the current sequence, the sequence will eventually reach $1$.
The conjecture remains open, but we will report on some numerical
results and conjectures in the case when $S$ consists of only $2$'s and $3$'s.

\vspace{0.8\baselineskip}



\section{The curling number conjecture}\label{Sec1}

Let $S$ be a finite nonempty sequence of integers.
By grouping adjacent terms, it is always possible
to write it as $S = X\,Y\,Y\,\ldots\,Y = X\,Y^k$, 
where $X$ and $Y$ are sequences of integers 
and $Y$ is nonempty ($X$ is allowed to be the empty sequence $\emptyset$).
There may be several ways to do this: choose the one
that maximizes the value of $k$: this
$k$ is the {\em curling number} of $S$,
denoted by $k(S)$.

For example, if $S = 0\,1\,2\,2\,1\,2\,2\,1\,2\,2$,
we could write it as $X\,Y^2$, where $X = 0\,1\,2\,2\,1\,2\,2\,1$
and $Y = 2$, or as $X\,Y^3$, where $X = 0$
and $Y = 1\,2\,2$. The latter representation is to be
preferred, since it has $k=3$, and as $k=4$
is impossible, the curling number of this $S$ is $3$.

Then we have

\begin{conj}\label{CNCconj}
The Curling Number Conjecture 
(\cite{GIJ}; also \cite{G4G7}, \cite{G4G8}).
If one starts with any initial sequence of integers $S$,
and extends it by repeatedly appending the curling number of 
the current sequence, the sequence will eventually reach $1$.
\end{conj}

In other words,
if $S_0 = S$ is any finite nonempty sequence of integers,
and we define 
\beql{Eq1}
S_{m+1} := S_m ~ k(S_m) \mbox{~~for~} m \ge 0 \,,
\eeq
then the conjecture is that for some $t \ge 0$ we will have $k(S_t) = 1$. 

For example, suppose we start with $S_0 = 2\,3\,2\,3$.
By taking $X = \emptyset$, $Y = 2\,3$, we have $S_0 = Y^2$,
so $k(S_0) = 2$, and we get $S_1 = 2\,3\,2\,3\,2$.
By taking $X = 2$, $Y = 3\,2$ we 
get $k(S_1) = 2$, $S_2 = 2\,3\,2\,3\,2\,2$.
By taking $X = 2\,3\,2\,3$, $Y = 2$ we 
get $k(S_2) = 2$, $S_3 = 2\,3\,2\,3\,2\,2\,2$.
Again taking $X = 2\,3\,2\,3$, $Y = 2$ we 
get $k(S_3) = 3$, $S_4 = 2\,3\,2\,3\,2\,2\,2\,3$.
Now, unfortunately, it is impossible to write 
$S_4 = X\,Y^k$ with $k>1$, so $k(S_4) = 1$,
$S_5 = 2\,3\,2\,3\,2\,2\,2\,3\,1$, and we have reached a $1$,
as predicted by the conjecture.
(If we continue the sequence from this point, it appears to join Gijswijt's
sequence, discussed in Section~\ref{Sec3}.)

Some of the proofs in \cite{GIJ} could be shortened and
the results extended if the conjecture were known to be true.
All the available evidence suggests that the conjecture {\em is} true,
but it has so far resisted all attempts to prove it.

\paragraph{Notation.}
We usually separate the terms of a sequence by small spaces.
$Y^k$ means $Y\,Y\, \ldots\,Y$, where $Y$ is repeated $k$ times.
$\emptyset$ denotes the empty sequence.
The curling number of $S$ is denoted by $k(S)$.
For a starting sequence of length $n$, 
$S_0 \, := \, s_1 \, s_2 \, \ldots \, s_n$,
where the $s_i$ are integers,
define
$
S_{m+1} := S_m ~ k(S_m) = s_1 \, \ldots \ s_{n+m+1}
$
for $m \ge 0$.
Suppose $t \ge 0$ is the smallest number such that
$k(S_t)=1$. Then we define 
$S^{(e)} := S_t = s_1 \, \ldots \ s_{n+t} \,.$
We call $S^{(e)}$ the {\em extended version} of $S_0$,
and $\tau(S_0) := t$ the {\em tail length} of $S_0$.
If $S_0$ never reaches 1, $S^{(e)} = S_{\infty}$ and $\tau(S_0) = \infty$.


\section{Sequences of $2$'s and $3$'s}\label{Sec2}

One way to approach the conjecture is to consider
the simplest nontrivial case, where the initial
sequence $S_0$ contains only $2$'s and $3$'s, and to
see how far such a sequence can extend using the
rule \eqn{Eq1} before reaching a $1$.
Perhaps if one is sufficiently clever,
one can invent a starting sequence that 
never reaches a 1, which would disprove the conjecture.
Of course it cannot reach a number greater that 3, either, for then
the following term will be a 1. So the sequence must be bounded 
between 2 and 3 for ever.
Unfortunately, even this apparently simple case 
has resisted our attempts to solve it.
Later in this section
we will mention some slight evidence that suggests 
the conjecture is true. 
First we report on our numerical experiments.

Let $\mu(n)$ denote the maximal length that can be achieved before
a $1$ appears, for any starting sequence $S_0$ consisting of $n$ $2$'s and $3$'s
(so $\mu(n) = n + \max_{S_0} ~ \tau(S_0)$, taken over all $S_0$ of length $n$).
If a 1 is never reached, we set $\mu(n) = \infty$.
The Curling Number Conjecture would imply $\mu(n) < \infty$
for all $n$. 

\begin{figure}[!h]
\centerline{\includegraphics[angle=270, width=5in]{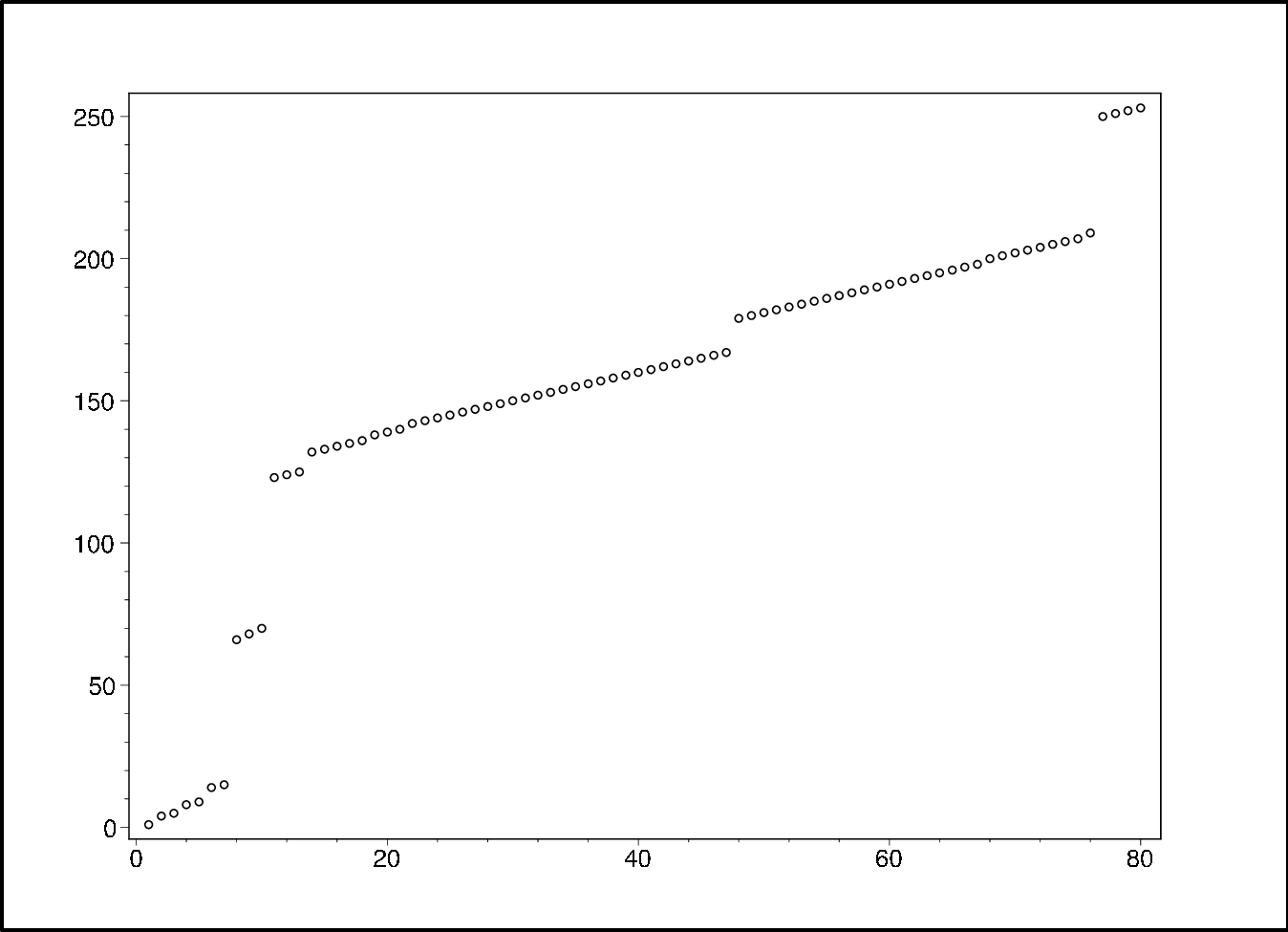}}
\caption{Scatter-plot of lower bounds on $\mu(n)$, the record for a starting sequence of $n$ $2$'s and $3$'s.
Entries for $n \le \LIMA$ are known to be exact; the other
entries are conjectured to be exact).}
\label{Fig1}
\end{figure}

By direct search, we have found $\mu(n)$ for all $n \le \LIMA$.
(The values for $n \le 30$ were given in \cite{GIJ}.)
The results are shown in Table \ref{Tabmu}
and Figure \ref{Fig1} together with lower bounds
(which we conjecture are in fact equal to $\mu(n)$)
for $\LIMB \le n \le \LIMC$.
The values of $\mu(n)$ also form sequence A094004 in \cite{OEIS}.

\begin{table}[htb]
$$
\begin{array}{|r|rrrrrrrrrrrr|}
\hline
     n & 1 & 2 & 3 & 4 & 5 & 6 & 7 & 8 & 9 & 10 & 11 & 12 \\
\mu(n) & 1 & 4 & 5 & 8 & 9 & 14 & 15 & 66 & 68 & 70 & 123 & 124 \\
\hline
     n & 13 & 14 & 15 & 16 & 17 & 18 & 19 & 20 & 21 & 22 & 23 & 24 \\
\mu(n) & 125 & 132 & 133 & 134 & 135 & 136 & 138 & 139 & 140 & 142 & 143 & 144 \\
\hline
     n & 25 & 26 & 27 & 28 & 29 & 30 & 31 & 32 & 33 & 34 & 35 & 36 \\
\mu(n) & 145 & 146 & 147 & 148 & 149 & 150 & 151 & 152 & 153 & 154 & 155 & 156 \\
\hline
     n & 37 & 38 & 39 & 40 & 41 & 42 & 43 & 44 & 45 & 46 & 47 & 48 \\
\mu(n) & 157 & 158 & 159 & 160 & 161 & 162 & 163 & 164 & 165 & 166 & 167 & 179 \\
\hline
     n &  49 &  50 &  51 &  52 &  53 &  54 &  55 &  56 &  57 &  58 &  59 &  60 \\
\mu(n) & 180 & 181 & 182 & 183 & 184^\dag & 185^\dag & 186^\dag & 187^\dag & 188^\dag & 189^\dag & 190^\dag & 191^\dag \\
\hline
     n &  61 &  62 &  63 &  64 &  65 &  66 &  67 &  68 &  69 &  70 &  71 &  72 \\
\mu(n) & 192^\dag & 193^\dag & 194^\dag & 195^\dag & 196^\dag & 197^\dag & 198^\dag & 200^\dag & 201^\dag & 202^\dag & 203^\dag & 204^\dag \\
\hline
     n &  73 &  74 &  75 &  76 &  77 &  78 & 79  & 80 &  &  &  &  \\
\mu(n) & 205^\dag & 206^\dag & 207^\dag & 209^\dag & 250^\dag & 251^\dag & 252^\dag  & 253^\dag  &  &  &  & \\
\hline
\end{array}
$$
\caption{$\mu(n)$, the record for extending a starting sequence of 
$n$ $2$'s and $3$'s before a 1 is reached. 
Entries marked with a dagger ($\dag$) are only lower bounds 
but are conjectured to be exact. }
\label{Tabmu}
\end{table}

In \cite{GIJ}, before we began computing $\mu(n)$,
we did not know how fast it would grow---would it
be a polynomial, exponential, or other function of $n$?
Even now we still do not know, since we have only limited
data. But up to $n = \LIMA$, and probably up to $n = \LIMC$,
$\mu(n)$ is a piecewise linear function of $n$.
There are occasional {\em jump points}, where
$\mu(n) > \mu(n-1)+1$, but in between jump points $\mu(n)$
simply increases by 1 when $n$ increases by $1$. 
Of course this piecewise linear behavior is not
incompatible with polynomial or exponential growth,
if the jump points are close enough together, but 
up to $n = \LIMC$ this seems not to be the case. 
There are long stretches of linear behavior.

The jump points are at 
$n = 1, 2, 4, 6, 8, 9, 10, 11, 14, 19, 22, 48$
and we believe the next three values are 68, 76 and 77
(entry A160766 in \cite{OEIS}).

From $n=2$ through $\LIMA$ (and probably through $n=\LIMC$)
the starting sequences $S_0$ which achieve $\mu(n)$
at the jump points are unique.
These especially good starting sequences are listed in Table \ref{TabGood1}.
For $2 \le n \le \LIMA$ (and probably
for $2 \le n \le \LIMC$)
these sequences $S_0$ have the following properties:
\begin{list}{--}{\setlength{\itemsep}{0.02in}}
\item
$S_0$ begins with 2,
\item
$S_0$ contains no nonempty subsequence of the form $W^4$,
\item
$S_0$ does not contain the subsequence 3,\,3 \,.
\end{list}

These are empirical observations. However, since they certainly hold
for the first $2^{\LIMB} -1$ choices for $S_0$,
we venture to make the following
\begin{conj}\label{Conj33}
If a starting sequence $S_0$ of length $n \ge 2$ achieves
$\mu(n)$ and $\mu(n) > \mu(n-1)+1$, then $S_0$ is unique
and has the above three properties.
\end{conj}

We can at least prove one thing about these especially
good starting sequences.
Let $S_0 \, := \, s_1 \, s_2 \, \ldots \, s_n$,
be any sequence of integers with extended version
$S^{(e)} := S_t = s_1 \, \ldots \ s_{n+t}$, where $k(S_t)=1$.
Call $S_0$ {\em weak} if each $S_r$ ($r=1, \ldots,t-1$)
can be written as $X \, Y^{s_{n+r}}$ with $X \ne \emptyset$.
In other words, $S_0$ is weak if $s_1$ is
not necessary for the computation
of the curling numbers $s_{n+1}, \ldots, s_{n+t}$.
This implies that $\tau(S_0) = \tau(\,s_2 \, \ldots \, s_n\,)$,
and establishes
\begin{lem}\label{LemWeak}
If a starting sequence $S_0$ of length $n \ge 2$ achieves
$\mu(n)$ and $\mu(n) > \mu(n-1)+1$, then $S_0$ is not weak.
\end{lem} 

\begin{table}[htb]
$$
\begin{array}{|l|l|}
\hline
n & \mbox{Starting sequence} \\
\hline 
1 & 2 \\
2 & 2\,2 \\
4 & 2\,3\,2\,3 \\
6 & 2\,2\,2\,3\,2\,2 \\
8 & 2\,3\,2\,2\,2\,3\,2\,3 \\
9 & 2\,2\,3\,2\,2\,2\,3\,2\,3 \\
10 & 2\,3\,2\,3\,2\,2\,2\,3\,2\,2 \\
11 & 2\,2\,3\,2\,3\,2\,2\,2\,3\,2\,2 \\
14 & 2\,2\,3\,2\,3\,2\,2\,2\,3\,2\,2\,3\,2\,3 \\
19 & 2\,2\,3\,2\,2\,3\,2\,3\,2\,2\,2\,3\,2\,2\,3\,2\,2\,3\,2 \\
22 & 2\,3\,2\,2\,3\,2\,2\,3\,2\,3\,2\,2\,2\,3\,2\,3\,2\,2\,3\,2\,2\,3 \\
48 & 2\,2\,3\,2\,2\,3\,2\,3\,2\,2\,2\,3\,2\,2\,2\,3\,2\,2\,3\,2\,2\,2\,3\,2\,2\,3\,2\,3\,2\,2\,2\,3\,2\,2\,2\,3\,2\,2\,3\,2\,2\,2\,3\,2\,2\,3\,2\,3 \\
\hline 
\end{array}
$$
\caption{Starting sequences of $n$ 2's and 3's for which $\mu(n) > \mu(n-1)+1$,
complete for $1 \le n \le \LIMA$.}
\label{TabGood1}
\end{table}

%

One further empirical observation is worth recording.
This concerns the starting sequences in between
jump points. Suppose $n_0$, $n_1$ are
consecutive jump points, so that
$$
\mu(n) ~=~ \mu(n-1) + 1 \quad \mbox{~for~} n_0 < n < n_1 \,,
$$
with $\mu(n) > \mu(n-1)+1$ at $n=n_0$ and $n_1$.
Then up to $n=\LIMA$ and conecturally up to $n=\LIMC$,
for $n_0 < n < n_1$, one
can obtain a starting sequence that achieves $\mu(n)$
by taking the starting sequence of length $n_0$
and prefixing it by a``neutral'' string of $n-n_0$
$2$'s and $3$'s
that do not get used in the computation of $\mu(n)$.
Although this is not surprising, we are unable to
prove that such neutral prefixes must always exist.
We return to this topic in Section \ref{SecRott}.

The large gaps between the jump points at 22 and 48 and between 48 and 68 are
especially noteworthy.
In particular, we have
\beql{Eq3}
\mu(n) = n + 120 \mbox{~for~} 22 \le n \le 47 \,,
\eeq
and, conjecturally,
\beql{Eq49}
\mu(n) = n + 131 \mbox{~for~} 48 \le n \le 67 \,.
\eeq

The data shown in Tables \ref{Tabmu}, \ref{TabGood1}
and Figure \ref{Fig1} for $n$ in the range $\LIMB$ to $\LIMC$
was obtained by computer search, assuming that the starting sequence
satisfied the three conditions listed in Conjecture \ref{Conj33},
but without assuming uniqueness. As it turned out, the best
starting sequences at the jump points were indeed unique.
The starting sequence we discovered at length 77 
produces (conjecturally) another large jump in $\mu(n)$.

We also made use of a number of more obvious shortcuts,
such as not considering a starting sequence
$s_1 \, \ldots \, s_n$ if
$k(s_1 \, \ldots \, s_{n-1}) = s_n$,
since we may assume that we have already considered all
starting sequences of length $n-1$.

The lemma cuts down the number of starting sequences of $2$'s and $3$'s
that must be considered.
Even if we simply exclude sequences that contain
four consecutive $2$'s or four consecutive $3$'s, 
the number of length $n$ (see entry A135491 in \cite{OEIS}) drops from $2^n$ 
to constant$\cdot \alpha^n$, where $\alpha = 1.839\ldots$.

Inspection of the best starting sequences in Table \ref{TabGood1}
suggests they must satisfy another condition,
which however we have not been able to prove: namely that they
do not contain two consecutive $3$'s. This is true for 
all the best starting sequences of lengths $n \le \LIMA$.

Making the assumption (as yet unjustified) that we need only consider
starting sequences with at most three consecutive $2$'s and with 
no pair of adjacent $3$'s, we were able extend the search to $n=78$.
This produced three further jumps, at $n=68$, $76$ and $77$,
establishing that
$\mu(n) \ge n+132$ for $68 \le n \le 75$,
$\mu(n) \ge n+133$ for $n=76$ and
$\mu(n) \ge n+173$ for $77 \le n \le \LIMC$.
The corresponding starting sequences are shown in Table \ref{TabGood2}.

\begin{table}[htb]
$$
\begin{array}{|l|l|}
\hline
n & \mbox{Starting sequence} \\
\hline
68 & 2\,2\,3\,2\,2\,3\,2\,2\,2\,3\,2\,2\,3\,2\,3\,2\,2\,2\,3\,2\,2\,2\,3\,2\,2\,3\,2\,2\,2\,3\,2\,2\,3\,2\,2\,2\,3\,2\,2\,3\,2\,3\,2\,2\,2\,3\,2\,2 \\
 & \,2\,3\,2\,2\,3\,2\,2\,2\,3\,2\,2\,3\,2\,2\,2\,3\,2\,2\,3\,2 \\
76 & 2\,3\,2\,2\,2\,3\,2\,2\,3\,2\,2\,2\,3\,2\,2\,3\,2\,3\,2\,2\,2\,3\,2\,2\,2\,3\,2\,2\,3\,2\,2\,3\,2\,2\,2\,3\,2\,2\,3\,2\,2\,2\,3\,2\,2\,3\,2\,3 \\
 & \,2\,2\,2\,3\,2\,2\,2\,3\,2\,2\,3\,2\,2\,3\,2\,2\,2\,3\,2\,2\,3\,2\,2\,2\,3\,2\,2\,3 \\
77 & 2\,2\,3\,2\,2\,2\,3\,2\,3\,2\,2\,2\,3\,2\,2\,2\,3\,2\,2\,3\,2\,2\,2\,3\,2\,2\,2\,3\,2\,3\,2\,2\,2\,3\,2\,2\,2\,3\,2\,2\,3\,2\,2\,2\,3\,2\,2\,2 \\
 & \,3\,2\,3\,2\,2\,2\,3\,2\,2\,2\,3\,2\,2\,3\,2\,3\,2\,2\,2\,3\,2\,2\,3\,2\,2\,2\,3\,2\,3 \\
\hline
\end{array}
$$
\caption{Starting sequences of $n$ 2's and 3's for which $\mu(n) > \mu(n-1)+1$,
conjectured to be complete for $\LIMB \le n \le \LIMC$.}
\label{TabGood2}
\end{table}

We have not succeeded in finding any algebraic constructions
for good starting sequences.  However,
one simple construction enables us to
obtain lower bounds on $\mu(n)$ for some larger values of $n$.
Let $S$ be a sequence of length $n$ that achieves $\mu(n)$,
and let $T$ be the sequence of length $\mu(n)$ that
it generates (up to just before the first $1$ appears).
Then in some cases the starting sequence $T \, S$
will extend to $T \, T \, 2$ and beyond before reaching a $1$.
For example, taking $S$ to be the length-$48$ 
sequence in Table \ref{TabGood1},
the sequence $T \, S$ has length $179+48=227$ and extends to 
a total length of $596$ before reaching a $1$,
showing that $\mu(227)  \ge 596$.

It would be nice to have some further exact values of $\mu(n)$,
even though they will require extensive computations.
Can the especially good starting sequences shown in Tables
\ref{TabGood1} and \ref{TabGood2}
(especially at lengths 22, 48 and 77) be generalized?  
What makes them so special?
And above all, what is the asymptotic behavior of $\mu(n)$?


\section{Gijswijt's sequence}\label{SecRott}

\section{Gijswijt's sequence}\label{Sec3}

If we simply start with $S_0 = 1$, and generate an infinite sequence
by continually appending the curling number of the current
sequence, as in \eqn{Eq1}, we obtain
$$
1,1,2,1,1,2,2,2,3,1,1,2,1,1,2,2,2,3,2,1,1,2,1,1,2,2,2,3,1,1,2,1,1, \ldots \,.
$$
This is {\em Gijswijt's sequence},
invented by Dion Gijswijt when he was a graduate student
at the University of Amsterdam,
and analyzed in \cite{GIJ}.
It is entry A090822 in \cite{OEIS}.

The first time a $4$ appears is at term $220$. One can calculate
quite a few million terms without finding a $5$ (as
the authors of \cite{GIJ} discovered!), 
but in \cite{GIJ} it was shown that a $5$ 
eventually appears for the first time at about term
$$
10^{10^{\Snj 23}} \, .
$$
Reference \cite{GIJ} also shows that the
sequence is in fact unbounded, and conjectures that the first time that
a number $m~\ge 6$ appears is at about term number
$$
2^{2^{\Snj 3^{\Snj 4^{\cdot^{\cdot^{\cdot^{\Snj {m-1}}}}}}}} \,,
$$
a tower of height $m-1$.

There is another question, also still open, which relates
Gijswijt's sequence to the discussion in the previous section.
If we start with an initial sequence $S$ of $2$'s and $3$'s,
extend it until we reach the first $1$,
say at the $(k+1)$st step, and then keep going,
it {\em appears} that the result is always simply 
the first $k$ terms of the extension of $S$,
followed by Gijswijt's sequence.  In other words, there is never
any interaction between the first $k$ terms of the extension of $S$
and an initial segment of Gijswijt's sequence when computing 
curling numbers after the $k$th step. 
This seems plausible, but we would not be surprised if 
there was a counterexample.
It would be nice to have this question
settled one way or the other.

One final remark: To avoid $1$'s in the sequence,
we might define $h(S) = \max \{ k(S), 2 \}$,
and replace the recurrence \eqn{Eq1} by
\beql{EqH}
S_{n+1} = S_n ~ h(S_n) \mbox{~~for~} n \ge 0 \,,
\eeq
If we start with $S_0 = 1$ and use the rule \eqn{EqH}
to extend it, 
the resulting sequence (A091787) is again unbounded,
and now it is possible to compute exactly when the first $5$ appears,
which is at step
$$
77709404388415370160829246932345692180\,.
$$
See \cite{GIJ} for further information.

\end{document}